\documentclass[runningheads, envcountsame, a4paper]{llncs}

\usepackage{amssymb}
\usepackage{graphicx}
\usepackage{enumerate}
\usepackage{amsmath,amsopn,amsfonts}
\usepackage[hyphens]{url}
\usepackage{microtype}

\usepackage{url}
\urldef{\mailsa}\path| berdinsky@gmail.com, 
 prohrakju@gmail.com | 

\newcommand{\keywords}[1]{\par\addvspace\baselineskip
\noindent\keywordname\enspace\ignorespaces#1}

\begin{document}

\mainmatter  

\title{Nonstandard Cayley Automatic Representations   
	   for  Fundamental Groups of Torus 
       Bundles  Over the Circle
   }

\titlerunning{Nonstandard Cayley Automatic Representations 
	for Fundamental Groups}

%
%
\author{Dmitry Berdinsky\inst{1,2}
 \and 
        Prohrak Kruengthomya\inst{1,2} 
} 

\authorrunning{D. Berdinsky and 
               P. Kruengthomya}


\institute{Department of Mathematics, 
  Faculty of Science, 
  Mahidol University,  
 Bangkok, Thailand  \and 
 Centre of Excellence in Mathematics, 
 Commission on Higher Education, Bangkok, Thailand  
 \\
 \mailsa\\
}

%
%

\toctitle{Nonstandard Cayley Automatic Representations
	of Fundamental Groups} 
\tocauthor{Dmitry~Berdinsky and Prohrak Kruengthomya}
\maketitle
\setcounter{footnote}{0}

\begin{abstract}

  We construct a new family of Cayley automatic representations 
  of semidirect products 
  $\mathbb{Z}^n \rtimes_A \mathbb{Z}$ 
  for which none of the projections of 
  the normal subgroup $\mathbb{Z}^n$ onto 
  each of its cyclic components is finite automaton recognizable.  
  For $n=2$ we describe a 
  family of matrices from $\mathrm{GL}(2,\mathbb{Z})$ 
  corresponding to these representations. 
  We are motivated by a 
  problem of characterization of 
  all possible Cayley automatic representations of 
  these groups. 
 


\end{abstract}

 \section{Introduction and Preliminaries}
 \label{introsec}

   Thurston and Epstein showed that 
   a fundamental group of a closed 
   $3$--manifold is automatic if and only if
   none of its prime factors is a closed 
   manifold modelled on nilgeometry 
   or solvgeometry \cite[Chapter~12]{Epsteinbook}.
   A fundamental group of a 
   closed manifold modelled on nilgeometry or 
   solvgeometry has a finite index subgroup 
   isomorphic to $\mathbb{Z}^2 \rtimes_A \mathbb{Z}$, 
   where $A$ is unipotent or Anosov, 
   respectively. These groups are not 
   automatic due to 
   \cite[Theorems~8.2.8~and~8.1.3]{Epsteinbook}.
   To include all fundamental groups of 
   closed $3$--manifolds, the class of automatic groups 
   had been extended by Bridson and Gilman
   \cite{BridsonGilman1996}, Baumslag, Shapiro and 
   Short \cite{BSS99}; see also autostackable groups
   proposed by Brittenham, Hermiller and Holt 
   \cite{BrittenhamHermillerHolt14}. 
   In this paper we use the concept of Cayley automatic groups,  
   extending the class of automatic groups,   
   proposed 
   by Kharlampovich, Khoussainov 
   and Miasnikov \cite{KKM11}. 
   
   All semidirect products of the form 
   $\mathbb{Z}^n \rtimes_A \mathbb{Z}$ 
   are Cayley automatic \cite[Proposition~13.5]{KKM11}.  
   These groups are the fundamental groups of 
   torus bundles over the circle and they play 
   important role in  group theory.
   Bridson and Gersten studied 
   the Dehn function for this family groups
   \cite{BridsonGersten96}.       
   In this paper we construct a new 
   family of Cayley automatic representations       
   for semidirect products 
   $\mathbb{Z}^n \rtimes_A \mathbb{Z}$. 
   These representations 
   demonstrate unforeseen behaviour violating 
   a basic property, to be explained below in this section, 
   known for representations described in
   \cite[Proposition~10.5]{KKM11}. They also 
   reveal an unexpected connection with Pell's equation. 
   The results of this paper are based on the 
   original construction of FA--presentation for 
   $\left(\mathbb{Z}^2,+\right)$ found by Nies and 
   Semukhin \cite{NiesSemukhin07}.  
   
   In general, we are interested in the 
   following question: 
   Given a Cayley automatic group, 
   is there any way to characterize 
   all of its
   Cayley automatic representations in terms of some 
   numerical characteristics or by any other means? 
   Despite the generality of the 
   notion of Cayley 
   automatic groups which retains only computational 
   mechanism of automatic 
   groups, it is possible to partly answer 
   this question for some Cayley automatic groups 
   in terms of a certain numerical characteristic
   which is intimately related to the Dehn function. 
   We discuss it in more details in the end of 
   this section. In the following few paragraphs we
   briefly recall the notion of 
   Cayley automatic groups and representations, 
   and a standard way to construct 
   such representations for semidirect
   products $\mathbb{Z}^n \rtimes_A \mathbb{Z}$.        
   
   Let  $\Sigma$ be a finite alphabet.    
   We denote by $\Sigma_\diamond$  the alpahbet $\Sigma \cup \{\diamond \}$,  where 
   $\diamond \notin \Sigma$ is called a padding 
   symbol. The convolution 
   $w_1 \otimes \dots \otimes w_m \in \Sigma_\diamond ^m$ of 
   strings $w_1,\dots, w_m \ \in \Sigma^*$ is the string		of length $\max\{|w_1|,\dots,|w_m|\}$ obtained 
   by placing $w_1, \dots, w_m$ one under another 
   and adding the padding symbol $\diamond$ at the 
   end of each string to make their lengths equal.   
   More formally,  
   the $k$th symbol
   of $w_1 \otimes \dots \otimes w_m$ is $(\sigma_1,\dots,\sigma_m )^\top$, 
   where $\sigma_i$, $i=1,\dots,m$ is the $k$th symbol of 
   $w_i$ if $k \leqslant |w_i|$ and
   $\diamond$ otherwise. 
   The convolution $\otimes R$  of a $m$--ary relation 
   $R \subseteq \Sigma^{*m} $ is defined as
   $\otimes R = \{w_1 \otimes \dots \otimes w_m \, | \,
   (w_1, \dots, w_m) \in R\}$.
   The relation $R$ is called  
   FA--recognizable  
   if $\otimes R$ is recognized by a finite automaton. 
  
   Let $\mathcal{A} = (A; R_1 ^{m_1}, \dots, 
   R_\ell ^{m_\ell}, f_1 ^{k_1}, \dots,f_r ^{k_r})$ be a 
   structure, where $A$ is the domain,  
   $R_i ^{m_i} \subseteq A^{m_i}, 
   i=1,\dots,\ell$ is a $m_i$--ary relation over $A$
   and $f_j ^{k_j}: A^{k_j} \rightarrow A$, $j=1,\dots,r$ is
   a $k_j$--ary operation on $A$. 
   Assume that there exist 
   a regular language  $L \subseteq \Sigma ^*$ 
   and a bijection 
   $\psi: L \rightarrow A$ such that 
   all relations 
   $\psi^{-1}(R_i ^{m_i}) = \{ (w_1,\dots, w_{m_i}) 
   \in \Sigma^{*m_i}\, |\, (\psi (w_1), \dots, 
   \psi (w_{m_i}))  \in R_i ^{m_i} \}$,  
   $i = 1, \dots, \ell$ and  
   $\psi^{-1} (\mathrm{Graph} (f_j)) = \{ (w_1,\dots, w_{k_j}, w_{k_j +1}) \in \Sigma ^{*(k_j+1)}\, | \, f_j (\psi(w_1),\dots,\psi(w_{k_j}))$ $=
     \psi(w_{k_j +1})\}$, $j=1,\dots,r$
   are FA--recognizable. In this case the structure 
   $\mathcal{A}$ is called FA--presentable 
   and the bijection $\psi : L \rightarrow A$ is 
   called FA--presentation 
   of $\mathcal{A}$ 
   \cite{KhoussainovNerode95,Blumensath99,BakhFrankSasha2007}.  
   For a recent survey of the theory of FA--presentable 
   structures we refer the reader to \cite{Stephan2015}. 
   A finitely generated group $G$ is called Cayley automatic 
   if the labelled directed Cayley graph 
   $\Gamma (G,S)$ is a FA--presentable structure 
   for some generating set $S \subseteq G$~\cite{KKM11}.
   Cayley automatic groups form a special class of 
   FA--presentable structures and   
   they naturally generalize 
   automatic groups 
   retaining its basic algorithmic properties. 
   We call a FA--presentation $\psi: L \rightarrow G$ 
   of $\Gamma(G,S)$ a Cayley automatic 
   representation of the group $G$.

   We recall that every element of 
   a group $\mathbb{Z}^n \rtimes_A \mathbb{Z}$, 
   where $A \in \mathrm{GL}(n,\mathbb{Z})$,  
   is given as a pair 
   $(b,h)$, where $b \in \mathbb{Z}$ and  
   $h \in \mathbb{Z}^n$. The group 
   multiplication is given by 
   $(b_1,h_1) \cdot (b_2,h_2) = 
    (b_1 + b_2, A^{b_2}h_1 + h_2)$. 
   The maps $b \mapsto (b,\bf{0}) $ and 
   $h \mapsto (0,h)$ give the natural embeddings
   of $\mathbb{Z}$ and $\mathbb{Z}^n$ into 
   $\mathbb{Z}^n \rtimes_A \mathbb{Z}$, respectively, 
   where $0 $ and $\bf{0}$ denote the identities 
   of the groups $\mathbb{Z}$ and $\mathbb{Z}^n$, 
   respectively.          
   Let $g_0 = (1,\bf{0})$ and 
   $g_i = (0,e_i)$, where
   $e_i = (0,\dots,0,\underset{i}{1},0,\dots,0)^t 
   \in \mathbb{Z}^n$. 
   The elements $g_0,g_1,\dots,g_n$ generate the 
   group $\mathbb{Z}^n \rtimes_A \mathbb{Z}$. 
   The right multiplication by $g_i, i=0,1,\dots,n$ is
   as follows: for a given $g = (b,h) \in 
   \mathbb{Z}^n \rtimes_A \mathbb{Z}$, 
   $gg_0 = (b+1,Ah)$ and $gg_i = (b,h+e_i)$.    
  
   Let $\psi_1: L_1 \rightarrow \mathbb{Z}$ 
   be a Cayley automatic representation of 
   $\mathbb{Z}$  and
   $\psi_2: L_2 \rightarrow \mathbb{Z}^n$ be 
   a Cayley automatic representations 
   of $\mathbb{Z}^n$ such that the automorphism  of $\mathbb{Z}^n$ given by the matrix $A$ is FA--recognizable. 
   Then, due to \cite[Theorem~10.3]{KKM11}, one gets a Cayley automatic representation 
   $\psi: L \rightarrow 
    \mathbb{Z}^n \rtimes_A \mathbb{Z}$ as follows: 
   $L = L_1 L_2$ 
    (we may assume that 
     $L_1 \subset \Sigma_1$, 
     $L_2 \subset \Sigma_2$ and 
     $\Sigma_1 \cap \Sigma_2 = 
     \varnothing$) and 
     for given $u \in L_1$ and 
     $v \in L_2$, 
     $\psi(uv) = (\psi_1 (u), \psi_2 (v))$.
     A standard way to construct 
     $\psi_2: L_2 \rightarrow \mathbb{Z}^n$ 
     is to take a FA--presentation 
     $\varphi: L_0 \rightarrow \mathbb{Z}$ of
     the structure $(\mathbb{Z},+)$, for example
     a binary representation, and define 
     $L_2$ as 
     $L_2 = \{w_1 \otimes \dots \otimes w_n \, 
     | \, w_i \in  L_0, i=1,\dots,n\}$ and 
     $\psi_2$ as 
     $\psi_2 (w_1 \otimes \dots \otimes w_n) =
      (\varphi (w_1),\dots,\varphi (w_n))$ for 
      every $w_1,\dots,w_n \in L_0$. 
     Clearly, for such a representation $\psi_2$  
     every automorphism of $\mathbb{Z}^n$
     is FA--recognizable. 
     Therefore,  $\psi_1$ and $\varphi$
     as above give a Cayley automatic 
     representation of 
     $\mathbb{Z}^n \rtimes_A \mathbb{Z}$. 
     We call such a representation 
     standard. Every standard 
     Cayley automatic representation 
     $\psi : L \rightarrow 
      \mathbb{Z}^n \rtimes_A \mathbb{Z}$ 
     satisfies the following basic properties: 
     \begin{enumerate}[a)] 
        \item The language 
         $L_{\mathbb{Z}^n} =\psi^{-1} (\mathbb{Z}^n)$ 
         of the strings representing 
         elements in the subgroup 
         $\mathbb{Z}^n \unlhd 
         \mathbb{Z}^n \rtimes_A \mathbb{Z}$ 
         is regular and the 
         relation $R_A=\{(u,v) \in 
         L_{\mathbb{Z}^n} \times L_{\mathbb{Z}^n}
         \,|\, A\psi(u)=\psi(v)\}$ is FA--recognizable. 
         \item For each projection 
         $p_i : \mathbb{Z}^n \rightarrow \mathbb{Z}^n$, $i=1,\dots,n$, on the $i$th component 
         given by 
         $p_i ((z_1,\dots,z_n))=(0,\dots,0,z_i,0,
         \dots,0)$ the relation 
         $P_i = \{(u,v) \in L_{\mathbb{Z}^n} 
         \times L_{\mathbb{Z}^n} \,|\, p_i \psi(u) = 
         \psi(v)\}$ is FA--recognizable.          
     \end{enumerate}	 
  
  In this paper we construct Cayley 
  automatic representations of groups 
  $\mathbb{Z}^n \rtimes_A \mathbb{Z}$ 
  for which the property a) holds but 
  the property b) does not hold
  -- in other words, these representations 
  are nonstandard. Namely, 
  in Section \ref{niessemukhinsection}
  we construct Cayley automatic representations of 
  $\mathbb{Z}^n$ for which every projection 
  $p_i: \mathbb{Z}^n \rightarrow \mathbb{Z}^n$, $i=1,\dots,n$ 
  is not FA--recognizable while some nontrivial 
  automorphisms $A \in \mathrm{GL}(n,\mathbb{Z})$
  are FA--recognizable. A family of these automorphisms 
  for the case $n=2$ is 
  described in Section \ref{automorphismssection1}. 
  Taking such a representation 
  as $\psi_2$ and an arbitrary Cayley automatic 
  representation 
  $\psi_1: L_1 \rightarrow \mathbb{Z}$ one 
  obtains a Cayley automatic representation 
  of $\mathbb{Z}^n \rtimes_A \mathbb{Z}$ as 
  described above. Clearly, for this representation
  the property a) holds    
  while the property b) 
  does not hold. In this paper we primarily focus on the case $n=2$  
  briefly discussing the case $n>2$.   Section 
  \ref{conclusionsec} concludes the paper.   
  
  Apart from the importance of 
  semidirect products $\mathbb{Z}^n \rtimes_A \mathbb{Z}$,  
  let us explain another   
  reason motivated us to study Cayley automatic representations of this family of groups
  violating 
  at least one of 
  the properties a) or b).  
  We first briefly recall some notation and results. 
  For a given f.g. group 
  $G$ with some finite set of generators 
  $A \subseteq G$,  we denote by $A^{-1}$ the set of inverses of 
  the elements of $A$ in $G$ and by $d_A$ the word metric
  in $G$ with respect to $A$. 
  We denote by 
  $\pi : \left(A \cup A^{-1}\right)^* \rightarrow G$ the 
  canonical map sending a word $w \in (A \cup A^{-1})^*$ 
  to the corresponding group element $\pi (w)$.  
  For the rest of the section we assume that $L \subseteq (A \cup A^{-1})^*$
  \footnote{We recall that every FA--presentable structure 
  has a FA--presentation over a binary alphabet 
  \cite{Blumensath99}. The alphabet $A \cup A^{-1}$ always
  has at least two symbols. The case of 
  FA--presentable structures over a unary alphabet 
  is special, see 
  \cite{Blumensath99,UnaryAutomaticGraphs01,UnaryAutomaticGraphs08}.}. 
  We denote by  $L^{\leqslant n}$ the language 
  $L^{\leqslant n} = 
  \{w \in  L \mid |w| \leqslant n \}$.      
  For a Cayley automatic representation 
  $\psi : L \rightarrow G$ we denote by 
  $h$ the function: $h(n) = \max \{d_A (\psi(w), \pi (w)) | 
  w \in L^{\leqslant n}\}$. 
 The function $h$ had been 
  introduced in \cite{measuring_closeness1}  as 
  a measure of deviation of Cayley automatic representation
  $\psi$ from $\pi$, i.e., from being automatic in the
  classical sense of Thurston.  
  For two nondecreasing functions 
  $h: [Q_1,+\infty) \rightarrow \mathbb{R}^+$ and 
  $f: [Q_2,+\infty) \rightarrow \mathbb{R}^+$, 
  where $[Q_1,+\infty),[Q_2,+\infty) 
  \subseteq \mathbb{N}$, we say that 
  $h \preceq f$ if there exist positive 
  integer constants $K,M$ and $N$
  such that for all $n \geqslant N$: 
  $h (n) \leqslant K f(Mn)$.  
  A f.g. group is said to be in 
  $\mathcal{B}_f$ if there exists a Cayley automatic 
  representation $\psi$ for which the function $h \preceq f$.   
  It was shown that the identity function $\mathfrak{i}(n)=n$ 
  is the sharp lower bound of the function $h$
  (in the sense of $\preceq$)   
  for all Cayley automatic 
  representations of the Baumslag--Solitar 
  groups $BS(p,q), 1 \leqslant p <q$ 
  \cite[Theorem~11]{measuring_closeness1} and 
  the wreath products $G \wr H$, if
  $H$ is virtually cyclic and $G$ is in the class 
  $\mathcal{B}_{\mathfrak i}$ \cite{BET19}.
  
  We recall that 
  the Heisenberg group 
  $\mathcal{H}_3 (\mathbb{Z})$ is isomorphic 
  to $\mathbb{Z}^2 \rtimes_T \mathbb{Z}$  
  for some lower triangular matrix $T$, see Remark  
  \ref{Heisenberg_remark1}. 
  The result of  \cite[Theorem~5.1]{eastwest19} 
  shows that if a Cayley automatic
  representation of the Heisenberg group 
  $\psi: L \rightarrow \mathcal{H}_3(\mathbb{Z})$ satisfies
  certain conditions, then the function 
  $h$ is bounded from below by the exponential function  
  $\mathfrak{e}(n) = \exp (n)$. 
  In particular,  
  for every Cayley automatic representation 
  $\psi: L \rightarrow \mathcal{H}_3 (\mathbb{Z})$ 
  satisfying the properties a) and b) 
  the function $h$ has the exponential lower 
  bound: $\mathfrak{e} \preceq h$.
  The lower bounds  
  for all possible Cayley automatic 
  representations of the  Heisenberg group 
  and the groups $\mathbb{Z}^2 \rtimes_A \mathbb{Z}$, 
  if $A \in \mathrm{GL} (2,\mathbb{Z})$ is a 
  matrix with two real eigenvalues not equal to $\pm 1$,  
  known to us are given by the 
  functions $\sqrt[3]{n}$ and $\mathfrak{i}$, 
  respectively, see \cite[Corollary~2.4]{eastwest19}.  
  However, it is not known whether or not these lower bounds 
  are sharp.  
  These observations motivated us to seek 
  nonstandard Cayley automatic representations 
  for a whole family of groups 
  $\mathbb{Z}^n \rtimes_A \mathbb{Z}$, 
  $A \in \mathrm{GL}(n,\mathbb{Z})$. 
  While we construct nonstandard
  representations for a large family 
  of groups $\mathbb{Z}^n \rtimes_A \mathbb{Z}$, see Theorem 
  \ref{classificationtheorem1} for the case $n=2$, 
  it does not contain nilpotent groups including 
  the Heisenberg group 
  $\mathcal{H}_3 (\mathbb{Z})$. 
  This leads us 
  to think that the case of nilpotent groups 
  is special. 
                  
 \section{Nies--Semukhin   
          FA--presentations of $\left(\mathbb{Z}^n,+\right)$} 
  \label{niessemukhinsection}

  Nies and Semukhin constructed 
  a FA--presentation of 
  $\left(\mathbb{Z}^2,+\right)$ for which no nontrivial 
  cyclic subgroup is FA--recognizable
  \cite[\S~6]{NiesSemukhin07}. 
  Let us briefly recall their 
  construction.   
  The group  $\mathbb{Z}^2$ is identified 
  with the additive group of the 
  quotient ring 
  $\mathbb{Z}[x]/\langle p_3 \rangle$, where 
  $p_3 (x) = x^2 + x - 3$  
  \footnote{In \cite[Remark~6.1]{NiesSemukhin07}
  it is said that one can use a polynomial 
  $x^2 + x - q$ for a prime $q \geqslant 3$.}.    
  A polynomial 
  $a_n x^n + \dots + a_0 \in \mathbb{Z}[x]$ 
  is called reduced if $|a_i| \leqslant 2$ for all
  $i=0,\dots,n$.  
  For given 
  $f,g \in \mathbb{Z}[x]$, 
  it is said that   
  $f \sim g$ if $p_3$ divides $f-g$. 
  In \cite[Proposition~6.2]{NiesSemukhin07}
  it is then shown that 
  every $f(x) \in \mathbb{Z}[x]$ 
  is equivalent to a reduced polynomial
  $\widetilde{f}(x)$.  
  Let $\Sigma=\{-2,-1,0,1,2\}$.     
  Each reduced polynomial
  $a_n x^n + \dots + a_0$ is
  represented by a string  
  $a_0 \dots a_n$ over the alphabet $\Sigma$. 
  Two strings 
  $u = a_0 \dots a_n$ and  
  $v= b_0 \dots b_m$ from $\Sigma^{*}$ 
  are said to be equivalent ($u \sim v$) 
  if  
  $a_n x^n + \dots + a_0 \sim
   b_m x^m + \dots + b_0$. 
  It is then shown that
  this equivalence relation 
  defined on $\Sigma^*$  
  is FA--recognizable. 
  Let $llex$ be the
  length--lexicographical 
  order on $\Sigma^*$ with respect to 
  the ordering $-2<-1<0<1<2$. 
  A regular domain for a presentation of 
  $\mathbb{Z}^2$ is defined as 
  $ \mathrm{Dom} = 
    \{ w \in \Sigma^* : 
    (\forall u < _{llex}w) 
    \, u \not\sim   w \}$. 
  Then a FA--recognizable relation 
  $R(x_1,x_2,x_3) \subset \Sigma^{*3}$ is 
  defined such that 
  for every pair $x_1,x_2 \in \Sigma^{*}$ 
  there exists a unique $x_3 \in \Sigma^*$ for 
  which $(x_1,x_2,x_3) \in R$ and 
  if $(x_1,x_2,x_3) \in R$, 
  then for the corresponding 
  polynomials $f_1,f_2$ and $f_3$:
  $f_1 + f_2 \sim f_3$. 
  It enables to define 
  a FA--recognizable relation  
  $\mathrm{Add}(x,y,z)$ on 
  $\mathrm{Dom}$ as follows: 
  $\mathrm{Add} = \{(x,y,z): 
   x,y,z \in \mathrm{Dom} \wedge 
   \exists w (R(x,y,w) \wedge 
   (w \sim z))\}$. 
  Clearly, the structure 
  $(\mathrm{Dom},\mathrm{Add})$ is 
  isomorphic to $(\mathbb{Z}^2,+)$.

  Now we notice that
  the Nies--Semukhin construction 
  can be generalized for a given 
  polynomial $t(x) = x^2 + p x - q 
  \in \mathbb{Z}[x]$ 
  for which $1+ |p| < |q|$.  
  Again, we identify 
  $\mathbb{Z}^2$ with the additive
  group of the quotient 
  ring $\mathbb{Z}[x]/\langle t \rangle$.  
  The inequality $1+ |p| < |q|$ implies that 
  $|q| \geqslant 2$.  
  We say that a polynomial
  $a_n x^n + \dots + a_0 \in \mathbb{Z}[x]$ 
  is reduced if $|a_i|<|q|$ for all $i =0,\dots,n$ 
  and two polynomials $f,g \in \mathbb{Z}[x]$ 
  are equivalent $f \sim g$
  if $t$ divides $f-g$. 
  For a given real $r$ we denote 
  by $[r]$ the integral part of 
  $r$: $[r]= \max \{m \in \mathbb{Z} \,|\, 
  m \leqslant r \} $ if $r \geqslant 0$ 
  and 
  $[r]= \min\{ m \in \mathbb{Z} \, 
  | \, m \geqslant r \}$ if $r<0$. 
  \begin{proposition} 
  \label{everypolequivtored_prop1}     
     Every polynomial 
     $f(x) \in \mathbb{Z}[x]$ is 
     equivalent to a reduced polynomial 
     $\widetilde{f}(x)$. 
  \end{proposition}  
  \begin{proof} 
     Let     
     $f(x)= a_n x^n + \dots + a_0$
     and      
     $k_0 = \left[\frac{a_0}{q}\right]$. 
     Since 
     $x^2 + px \sim q$, 
     $f(x) \sim f_1(x) = b_n x^n + \dots + b_0$, 
     where 
     $b_0 = a_0 - 
     k_0 q$, 
     $b_1 = a_1 + k_0 p$, 
     $b_2 = a_2 + k_0$ and 
     $b_i = a_i$ for $i>2$.
     If $|a_0|<|q|$, then $f_1 (x) =f_0 (x)$. 
     Otherwise, we get that $\sum_{i=0}^n |a_i| >
     \sum_{i=0}^n |b_i|$.     
     Let $k_1=\left[\frac{b_1}{q}\right]$.    
     Since $x^3 + px^2 \sim qx$, 
     $f_1(x) \sim f_2(x)=c_n x^n + \dots +c_0$, 
     where $c_0 = b_0$, 
     $c_1=b_1 - k_1 q$, 
     $c_2 = b_2 + k_1 p$, 
     $c_3 = b_3 + k_1$ and 
     $c_i = b_i$ for $i>3$.  
     If $|b_1| < |q|$, then $f_2 (x) = f_1 (x)$. 
     Otherwise, we get that
     $\sum_{i=0}^n |b_i| >
     \sum_{i=0}^n |c_i|$. 
     We have: $|c_0|=|b_0|<|q|$ and 
     $|c_1|<|q|$. 
      If we continue in this way, the 
      process will terminate after 
      a finite number of iterations 
      producing 
      a reduced polynomial 
      $\widetilde{f}(x)$ at the last  
      iteration.\hfill\qed
  \end{proof}  
 \begin{remark}
 	It can be seen that if the inequality $1+|p|<|q|$ 
 	is not satisfied, then the procedure 
 	described in Proposition \ref{everypolequivtored_prop1}
 	fails to produce a reduced polynomial 
 	for some input polynomials $f(x)$.  
 	For example, let $t(x)=x^2 + 2x -3$ and 
 	$f(x) = 2 x + 6$. Applying the procedure
 	from Proposition \ref{everypolequivtored_prop1} 
 	one gets an infinite sequence of 
 	polynomials $f_i (x) = 2 x^{i+1} + 6 x^i$
 	which never terminates.  
 \end{remark} 	
  
  Let  
  $\Sigma_q = \{ -(|q|-1),  
  \dots, |q|-1 \}$. 
  We represent a reduced polynomial 
  $a_n x^n + \dots + a_0$ by a string 
  $a_0 \dots a_n$ over the alphabet 
  $\Sigma_q$. 
  Similarly, we say that two 
  strings $a_0 \dots a_n$ and 
  $b_0 \dots b_m$ over $\Sigma_q$ 
  are equivalent if the polynomials 
  $a_n x_n + \dots + a_0$ and  
  $b_m x^m + \dots + b_0$ are equivalent.  
  An algorithm checking whether 
  two given reduced polynomials 
  $f(x) = a_n x^n + \dots + a_0$ and 
  $g(x) = b_m x^m + \dots + b_0$ are 
  equivalent 
  is the same, up to minor changes,
  as it is 
  described by Nies and Semukhin 
  for the case $t(x)= x^2 + x -3$,   
  see \cite[\S~6]{NiesSemukhin07}. 
  We first check if 
  $q$ divides $a_0 - b_0$; if not, 
  $f \not\sim g$.   
  We remember two carries  
  $r_0 = p \frac{a_0-b_0}{q} $  and 
  $r_1 = \frac{a_0-b_0}{q}$, and  
  then 
  verify whether $q$ divides $r_0+a_1-b_1$;
  if not, $f \not\sim g$. Otherwise, 
  we update the carries: 
  $r_0 \rightarrow r_1 + 
  p\frac{r_0+a_1-b_1}{q}$ and 
  $r_1 \rightarrow \frac{r_0+a_1-b_1}{q}$, 
  and then verify whether 
  $q$ divides $r_0 + a_2 - b_2$. 
  Proceeding in this way we  
  check if $f \sim g$ or not.    
  Initially, $|r_1| \leqslant 1 
   \leqslant |q|-1$ and 
  $|r_0| \leqslant |p| < (|q|-1)^2$. 
  Since $q$ divides $r_0 + a_i - b_i$ 
  at every step of our process unless 
  $f \not\sim  g$, we can change the formulas 
  for updating carries as follows: 
  $r_0 \rightarrow r_1 + 
  p \left[ \frac{r_0+a_i-b_i}{q} 
  \right]$ and  
  $r_1 \rightarrow 
  \left[\frac{r_0 + a_i-b_i}{q}\right]$. 
  Now, if $|r_1| \leqslant |q|-1$ and 
  $|r_0| \leqslant (|q|-1)^2$, then
  $ 
  \left|\left[\frac{r_0 + a_i - b_i}{q}\right]
  \right| 
  \leqslant 
  \left[\frac{(|q|-1)^2 +2(|q|-1)}{|q|}\right]
  = |q|-1$ 
  and  
  $
  \left|r_1 + 
  p \left[ \frac{r_0+a_i-b_i}{q} 
  \right]\right| \leqslant  (|q|-1) + 
  |p|\left| 
  \left[ \frac{r_0+a_i-b_i}{q} \right]
  \right|
  \leqslant (|q|-1)+(|q|-2)(|q|-1)= 
  (|q|-1)^2 
  $. This shows that $|r_1|$ and $|r_0|$ 
  are always bounded by $|q|-1$ and 
  $(|q|-1)^2$. 
  This algorithm requires only a
  finite amount of memory, so the 
  equivalence relation $\sim$ is 
  FA--recognizable. 
  
  Similarly, one can construct a 
  FA--recognizable relation 
  $R(u,v,w) \subset \Sigma_q ^*$ 
  such that for every pair 
  $(u,v) \in \Sigma_q ^*$ there 
  exists a unique $w \in \Sigma_q ^*$ 
  for which $(u,v,w) \in  R$  
  and if $(u,v,w) \in R$ then 
  for the corresponding polynomials 
  $f_u,f_v$ and $f_w$: $f_u + f_v \sim f_w$.
  Again, the construction of such a relation
  $R$ is 
  the same, up to minor changes, 
  as it is described by Nies and Semukhin
  for the case $t(x)=x^2 + x - 3$. 
  Let  $u = a_0 \dots  a_n$ and 
  $v = b_0 \dots b_m$. Then a string 
  $w= c_0 \dots c_k$ for which 
  $(u,v,w) \in R$ is obtained as follows. 
  Let  $c_0$ be an integer such that 
  $|c_0| < |q|-1$, $c_0$ has the same sign 
  as $a_0 + b_0$  and 
  $c_0 \equiv a_0 + b_0
  \left(\mathrm{mod} \,q \right)$. We   
  remember two carries 
  $r_0 = p \left[\frac{a_0 + b_0}{q}\right]$ 
  and $r_1 = \left[\frac{a_0 + b_0}{q}\right]$. 
  We put $c_1$ to be an integer such that
  $|c_1| \leqslant |q|-1$, 
  $c_1$ has the same sign as 
  $r_0 + a_1 + b_1$ and 
  $c_1 \equiv r_0 + a_1 + b_1  
  \left(\mathrm{mod} \,q \right)$,  
  and update the carries as 
  $r_0 \rightarrow r_1 + 
  p \left[\frac{r_0 + a_1 + b_1}{q}\right]$ 
  and 
  $r_1 \rightarrow 
  \left[\frac{r_0 + a_1 + b_1}{q}\right]$. 
  This process is continued until the string 
  $w$ is generated. The formulas for 
  updating carries are 
  $r_0 \rightarrow 
   r_1 + p 
   \left[\frac{r_0 + a_i + b_i}{q}\right]$ 
  and $r_1 \rightarrow 
  \left[\frac{r_0 + a_i + b_i}{q}\right]
  $. The proof that $|r_1|$ and $|r_0|$ 
  are bounded by $(|q|-1)$ and $(|q|-1)^2$, 
  respectively, is the same as in the 
  paragraph above, so the relation $R$ is 
  FA--recognizable. 
  
  Fixing the ordering 
  $-(|q|-1)< \dots < (|q|-1)$ on 
  $\Sigma_q$, the domain 
  $\mathrm{Dom}$ and the relation 
  $\mathrm{Add}$ are then defined 
  in exactly the same way 
  as by Nies and Semuhkhin, see the 
  first paragraph of this section. 
  So, for every pair of integers 
  $p$ and $q$, for which $1 + |p| < |q|$,
  we obtain a regular domain 
  $\mathrm{Dom}_{p,q}$ and a FA--recognizable 
  relation $\mathrm{Add}_{p,q}$ 
  for which $(\mathrm{Dom}_{p,q}, \mathrm{Add}_{p,q})$ 
  is isomorphic to $(\mathbb{Z}^2,+)$. 
  For given $p$ and $q$ satisfying the inequality  
  $1+|p|<|q|$, we denote by 
  $\psi_{p,q}: 
  \mathrm{Dom}_{p,q} \rightarrow \mathbb{Z}^2$ 
  the representation of $(\mathbb{Z}^2,+)$ 
  described above. 
  Let $g  \in \mathbb{Z}[x]$ be 
  some fixed polynomial. Clearly, 
  if $f_1 \sim f_2$, then 
  $f_1 g \sim f_2 g$. 
  Therefore, multiplication by $g$ 
  induces a map from 
  $\mathbb{Z}[x]/ \langle t \rangle$ 
  to $\mathbb{Z}[x]/ \langle t \rangle$ 
  which sends  
  an equivalence class $[f]_{_\sim}$ 
  to the equivalence class 
  $[fg]_{_\sim}$. 
  So, by Proposition 
  \ref{everypolequivtored_prop1}, 
  multiplication by $g$ induces 
  a map $\varphi_g : \mathrm{Dom}_{p,q} 
  \rightarrow \mathrm{Dom}_{p,q}$.   
  \begin{proposition}
  \label{multiplication_by_g_FArec_prop1}   
     For every representation 
     $\psi_{p,q}$ the function  
     $\varphi_g : 
     \mathrm{Dom}_{p,q} \rightarrow 
     \mathrm{Dom}_{p,q}$ is 
     FA--recognizable. 
  \end{proposition} 
  \begin{proof} 
     Since the equivalence relation $\sim$ and 
     $\mathrm{Add}$ are FA--recognizable, 
     it is enough only 
     to show that multiplication 
     by a monomial $x$ is FA--recognizable. 
     It is true because for a string  
     $u = a_0 \dots a_n \in \mathrm{Dom}_{p,q}$ 
     the string 
     $\varphi_x (u)$ is equivalent 
     to the shifted string 
     $0 a_0 \dots a_n$. 
     Clearly, such shifting of strings is  
     FA--recognizable.\hfill\qed    
  \end{proof}
 
  Nies and Semukhin showed that 
  every nontrivial cyclic subgroup 
  $\langle z \rangle$ of 
  $\mathbb{Z}^2$ is not 
  FA--recognizable for the 
  representation $\psi_{1,3}$ 
  \cite[\S~6]{NiesSemukhin07}.     
  We will show that each of the two cyclic 
  components of $\mathbb{Z}^2$ is not 
  FA--recognizable for every representation
  $\psi_{p,q}$, if $\gcd(p,q)=1$.   
  Let $\xi = [1]_{_\sim}$, 
  where $1$ is the 
  polynomial $f(x)=1$; also, $\xi$ corresponds 
  to the single--letter string 
  $1 \in \mathrm{Dom}_{p,q}$: 
  $\psi_{p,q}(1)=\xi$. 
  Let us show that  
  the cyclic subgroup generated by $\xi$ is not 
  FA--recognizable with respect to 
  $\psi_{p,q}$, if $\gcd(p,q)=1$.
  We will use arguments analogous 
  to the ones in \cite[\S~6]{NiesSemukhin07} 
  with relevant modifications. 
  It is straightforward that
  \cite[Lemma~6.3]{NiesSemukhin07}
  claiming that 
  for given two equivalent reduced polynomials 
  $f(x)$ and $g(x)$,  
  $x^k | f$  implies 
  $x^k | g$, holds valid. 
  It is said that $f(x) \in \mathbb{Z}[x]$  
  starts with $k$ zeros in reduced 
  form if there exists a reduced 
  polynomial $g(x)$ for which  $f \sim g$ and  
  $x^k | g(x)$: in this case the string 
  representing $g(x)$ starts with  
  $k$ zeros. 
  For a given $k > 0$, the polynomial 
  $q^k$ starts with at least $k$ zeros 
  in reduced form because  
  $q^k \sim x^k(x+p)^k$. 
  
  Assume now that $L_\xi = 
  \psi_{p,q} ^{-1} (\langle \xi \rangle)$ is regular and
  recognized by a finite automaton with $k_0$ 
  states. 
  The string 
  $\psi_{p,q}^{-1}
   ([q^{k_0}]_{_\sim}) \in L_\xi$ 
  starts with at least $k_0$ zeros, 
  i.e., $\psi_{p,q}^{-1}
   ([q^{k_0}]_{_\sim}) = 0^ku$ for 
   $k \geqslant k_0$ and some 
   $u \in \Sigma_q^*$, which does not
   have $0$ as the first symbol.    
  By pumping lemma, 
  there exist $k_1,k_2$ and 
  $0 < d \leqslant k_0$,
  for which $k_1+d+k_2=k$, 
  such that $s_i = 0^{k_1 + d i + k_2} u 
  \in L_\xi$ for all $i\geqslant 0$. 
  Since $s_i \in L_\xi$, we have a 
  sequence of integers 
  $n_i$, $i \geqslant 0$ for which 
  $\psi_{p,q}(s_i) = [n_i]_{_\sim}$, 
  so $n_i$ starts with $k_1 + d i + k_2$ 
  zeros in reduced form. 
  For a given integer $n$, if it starts 
  with at least one zero in reduced 
  form, then $q\,|\,n$: 
  it is because  
  $n = q \ell + r$ 
  for some $\ell$ and 
  $r \in \{0,\dots,|q|-1\}$, 
  so if $r \neq 0$ then 
  $n \sim x(x+p) \ell + r$ starts 
  with no zeros in reduced form.  
  \begin{proposition}     
  \label{gcdpqequalone_prop}     
     Assume that $\gcd(p,q)=1$. 
     If $n = q \ell$ starts with $m>0$ zeros in 
     reduced form, then $\ell$ starts with 
     $m-1$ zeros in reduced form. 
  \end{proposition}  
  \begin{proof}   
     Let $f(x) = 
     x^i (b_j x^{j-i} + \dots + b_i)$ 
     be a reduced 
     polynomial equivalent to $\ell$, 
     where $b_i \neq 0$.       
     We have $n=q \ell \sim x^{i+1} (x+p)
     (b_j x^{j-i} + \dots + b_i)$. 
     Since $\gcd(p,q)=1$ and 
     $|b_i| < |q|$, 
     $q  \not| \,\, pb_i$. Therefore, 
     $n$ starts with $i+1$ zeros in 
     reduced form, so $i = m -1$. 
     Therefore, $\ell$ starts with 
     $m-1$ zeros in reduced form.\hfill\qed 
  \end{proof}

  Thus, if $\gcd (p,q)=1$, by Proposition 
  \ref{gcdpqequalone_prop}, we obtain that 
  $q^{k_1 + d i + k_2} | n_i$, 
  so $n_i = q^{k_1 + d i + k_2} m_i$ 
  for some nonzero integer $m_i$.   
  Let $\alpha$ and $\beta$ be the 
  roots of the polynomial 
  $t(x) = x^2 + px - q$. 
  We have $\alpha \beta = -q$, 
  so $|\alpha \beta| = |q|$. 
  Therefore, either $|\alpha|$ or 
  $|\beta|$ must be less 
  or equal than $\sqrt{|q|}$. 
  So, let us assume that 
  $|\alpha| \leqslant \sqrt{|q|}$.  
  For every two equivalent polynomials 
  $f \sim g$: $f(\alpha) = g(\alpha)$. 
  Let $f_i$ be the reduced polynomials 
  corresponding to the strings $s_i$. 
  If $|\alpha| > 1$, 
  then $|f_i(\alpha)|$ is bounded from 
  above by 
  $(|q|-1)|u||\alpha|^{|s_i|-1}$, where 
  $|s_i|= k_1+d i+k_2+|u|$ is the 
  length of the string $s_i$;  
  it is because there are only at most 
  $|u|$ nonzero coefficients of the 
  polynomial $f_i$ and 
  the absolute value of each of which 
  is less than or equal to $|q|-1$. 
  Therefore, 
  $|f_i(\alpha)|\leqslant C_1 |\alpha|^{di}$,  
  where $C_1 = (|q|-1)|u|
         |\alpha|^{k_1 + k_2 + |u|-1}$. 
  If $|\alpha| \leqslant 1$, then 
  $|f_i (\alpha)| \leqslant C_2$, 
  where $C_2 = (|q|-1)|u|$. 
  In both cases we obtain that  
  $|f_i(\alpha)| \leqslant C \sqrt{|q|}^{di}$ 
  for some constant $C$. 
  On the other hand, since $f_i \sim n_i$, 
  $f_i(\alpha)= n_i = 
   q^{k_1 + di +k_2} m_i$. 
 Therefore, 
 $|f_i(\alpha)| =  |q|^{k_1+di+k_2}|m_i|
 \geqslant |q|^{di}$. 
 Thus, we obtain that 
 $|q|^{di} \leqslant C \sqrt{|q|}^{di}$
 for all $i \geqslant 0$, which 
 apparently leads to a contradiction
 since $|q|>1$. Thus, $L_\xi$ is not regular. 
 
 Let $\eta = [x]_\sim$, where $x$ is 
 the polynomial $f(x)=x$; also, $\eta$  
 corresponds to the string
 $01 \in \mathrm{Dom}_{p,q}$: 
 $\psi_{p,q}(01)=\eta$. 
 Clearly, $\mathbb{Z}^2$ is the direct sum 
 of its cyclic subgroups $\langle \xi \rangle$
 and $\langle \eta \rangle$. 
 Let $L_\eta = \psi_{p,q}^{-1} (\langle \eta \rangle)$.
 We notice that $L_\xi = \{w \in \mathrm{Dom}_{p,q}\,|\, 
 \varphi_x (w) \in L_\eta\}$. 
 The inclusion 
 $L_\xi \subseteq \{w \in \mathrm{Dom}_{p,q}\,|\, 
 \varphi_x (w) \in L_\eta\}$ is straightforward. 
 For the inclusion $\{w \in \mathrm{Dom}_{p,q}\,|\, 
 \varphi_x (w) \in L_\eta\} \subseteq L_\xi$ it is enough 
 to notice that 
 if $\psi_{p,q} (w) = [sx + r]_\sim$, then 
 $\varphi_x(w) =  [x(sx+r)]_\sim  = [s (-px + q) + rx]_\sim 
 = [(r-sp)x + sq]_\sim$ which is equal to $[kx]_\sim$ for 
 some $k \in \mathbb{Z}$ only if $sq = 0$.   
 The map $\varphi_x: 
 \mathrm{Dom}_{p,q} \rightarrow \mathrm{Dom}_{p,q}$ 
 is FA--recognizable, by Proposition 
 \ref{multiplication_by_g_FArec_prop1}. So, 
 the regularity of $L_\eta$ implies the regularity 
 of $L_\xi$. Therefore, $L_\eta$ is not regular. 
 Clearly, the fact that $L_\xi$ and $L_\eta$ are 
 not regular implies that the projections 
 of $\mathbb{Z}^2$ onto its cyclic components
 $\langle \xi \rangle$ and $\langle \eta \rangle$
 are not FA--recognizable.  
 Let us summarize the results we obtained 
 in the following theorem.
 
\begin{theorem}
 \label{niessem_thm1}    
    
 For every pair of integers $p$
 and $q$ for which $1+ |p| < |q|$ 
 the map $\psi_{p,q}: \mathrm{Dom}_{p,q}
 \rightarrow \mathbb{Z}^2$ gives 
 a FA--presentation of $(\mathbb{Z}^2,+)$.  
 Moreover, if $\gcd(p,q)=1$, then 
 none of the two cyclic components of 
 $\mathbb{Z}^2$ and the projections onto 
 theses components is FA--recognizable  
 with respect to $\psi_{p,q}$.   
\end{theorem}
\begin{remark}
\label{mult_by_constant_remark}	
	Let $z(x) = ax + b$ be a polynomial in $\mathbb{Z}[x]$  
	and $\zeta = [z]_\sim$ be the corresponding element in 
	$\mathbb{Z}^2$. For a given integer $m>1$ let  
	$\delta  = m \zeta \in \mathbb{Z}^2$. 
	We denote by $\varphi_m$ the map $\varphi_g$ for 
	the constant polynomial $g(x)=m$. 
	Let $L_\zeta  = \psi_{p,q}^{-1} 
	     (\langle \zeta \rangle)$ and 
	    $L_\delta = \psi_{p,q}^{-1} 
	     (\langle \delta \rangle)$.  
	Then  we have: $L_\zeta = \{ w \in 
	\mathrm{Dom}_{p,q} \, | \, 
	\varphi_{m} (w) \in L_\delta\}$. 
	The inclusion 
	$L_\zeta \subseteq 
	\{ w \in \mathrm{Dom}_{p,q} \, | \,
	 \varphi_m (w) \in L_\delta\}$ is straightforward. 
	In order to prove the inclusion 
	$\{w \in \mathrm{Dom}_{p,q}\, | \, 
	 \varphi_m (w) \in L_\delta \} \subseteq L_\zeta$ 
	we notice that if $\psi_{p,q}(w) = [sx+r]_\sim$, 
	then $\varphi_m  (w) = [m (sx + r)]_\sim$ which is equal 
	to $k \delta = [kmz]_\sim= [kmax + kmb]_\sim$ for some 
	$k \in \mathbb{Z}$ 
	iff $kma = ms$ and $kmb = mr$. 
	Clearly, this holds iff 
	$s= ka $ and $r = kb$, so 
	$\psi_{p,q}(w)= [kax+kb]_\sim = 
	 k \zeta$ which implies that $w \in L_\zeta$.
	Therefore, if $L_\zeta$ is not regular, then
	$L_\delta$ is not regular. 
	In particular, if $\gcd (p,q)=1$, then 
	none of the cyclic subgroups $\langle m \xi \rangle$ 
	and $\langle m \eta \rangle$ is FA--recognizable 
	with respect to $\psi_{p,q}$ 
	for every nonzero integer $m$.         
\end{remark}	
\begin{remark}    
\label{all_cyclic_subgroups_notFA} 	
    In order to guarantee that all 
    nontrivial cyclic 
    subgroups of $\mathbb{Z}^2$ are not 
    FA--recognizable with respect to 
    $\psi_{p,q}$, one should additionally 
    require that the polynomial 
    $t(x) = x^2 + px -q$ is 
    irreducible in $\mathbb{Z}[x]$. 
    Let $\gamma = [g]_\sim$ for 
    some $g \in \mathbb{Z}[x]$,  $g \not\sim 0$, 
    and $L_\gamma = \psi_{p,q}^{-1}(\langle \gamma
    \rangle)$.    
    We have: $L_\xi = \{ w \in \mathrm{Dom}_{p,q} \,
    |\,\varphi_g(w) \in  L_\gamma \}$.  
    The inclusion 
    $
     L_\xi \subseteq \{ w \in \mathrm{Dom}_{p,q} \, | 
     \, \varphi_g (w) \in L_\gamma \}
    $ is again straightforward. 
    In order to prove the inclusion 
    $\{w \in \mathrm{Dom}_{p,q} \, 
    |\, \varphi_g (w) \in L_\gamma \} \subseteq L_\xi$ 
    we notice that if $\psi_{p,q}(w) = [sx+r]_\sim$, 
    then $\varphi_g (w) = [g (sx+r)]_\sim$ which is 
    equal to $[g k]_\sim$ for some $k \in \mathbb{Z}$ 
    iff the polynomial $t$ divides 
    $g (sx + r -k)$. Since $t$ is irreducible 
    and $t$ does not divide $g$, then  $s=0$ and $r = k$.
    Therefore, by Proposition 
    \ref{multiplication_by_g_FArec_prop1}, 
    if $L_\gamma$ is regular, then $L_\xi$ is regular. 
    So, $L_\gamma$ is not regular.         	
 \end{remark}
 \begin{remark} 
 \label{degenerateendomorphisms_remark1} 	
 	Moreover, if $t$ is irreducible in $\mathbb{Z}[x]$,
 	every nonzero endomorphism of $\mathbb{Z}^2$ 
 	with nontrivial kernel is not FA--recognizable.
 	This immediately follows from the observation 
 	that the image of a such endomorphism is a 
 	cyclic subgroup of $\mathbb{Z}^2$ 
 	which is not FA--recognizable 
 	with respect to $\psi_{p,q}$ by  
 	Remark \ref{all_cyclic_subgroups_notFA}.   
 \end{remark}	

 Now, let $n>2$ and $t(x) = x^{n}+p_{n-1}x^{n-1}+\dots+p_{1}x-q$ be a polynomial with integer coefficients for which 
 $1+|p_{n-1}|+\dots+|p_{1}|<|q|$. We identify the group $\mathbb{Z}^n$ with the additive group of the ring 
 $\mathbb{Z}[x]/\langle t \rangle$. 
 We denote by $\overline{p}$ a tuple 
 $\overline{p} = \langle p_1, \dots, p_{n-1}  \rangle$. 
 Clearly, one gets a FA--presentation 
 $\psi_{\overline{p},q}: \mathrm{Dom}_{\overline{p},q} \rightarrow \mathbb{Z}^n$  of $\left(\mathbb{Z}^n,+\right)$ 
 in exactly the same way as it is described for the case $n=2$. 
 It can be seen that all arguments presented in this section 
 hold valid up to the following minor modifications. 
 For an algorithm recognizing the equivalence $\sim$, 
 one should use $n$ carries $r_0,r_1,\dots,r_{n-1}$ 
 updated as follows: 
 $r_{0}\rightarrow r_{1}+p_{1}[\frac{r_{0}+a_{i}-b_{i}}{q}], r_{1}\rightarrow r_{2}+p_{2}[\frac{r_{0}+a_{i}-b_{i}}{q}],...,r_{n-2}
 \rightarrow r_{n-1}+p_{n-1}[\frac{r_{0}+a_{i}-b_{i}}{q}],
 r_{n-1} \rightarrow[\frac{r_{0}+a_{i}-b_{i}}{q}]$. 
 Let us verify the inequalities: $r_0 \leqslant (|q|-1)^2, 
 r_1 \leqslant (|q|-1)(1+|p_{n-1}|+|p_{n-2}|+\dots+|p_{2}|), 
 \dots,|r_{n-2}|\leqslant(|q|-1)(1+|p_{n-1}|)$ and 
 $|r_{n-1}|\leqslant |q|-1$. 
 Initially these inequalities 
 are satisfied.  
 
 Suppose now that they hold for a current 
 iteration of the algorithm.   
 Since $|r_0| \leqslant (|q|-1)^2$ for the current 
 iteration, then 
 $\left| \left[\frac{r_0 + a_i - b_i}{q}\right] \right| 
  \leqslant |q|-1$ which implies that 
 $|r_{n-1}| \leqslant |q|-1$ for the next iteration of the 
 algorithm. 
 Since we assumed that $|r_{n-1}| \leqslant |q|-1$ for the 
 current iteration, then 
 $r_{n-1} + p_{n-1} \left[\frac{r_0 + a_i - b_i}{q}\right]
 \leqslant (|q|-1) (1 + |p_{n-1}|)$ which implies that
 $|r_{n-2}| \leqslant (|q|-1)(1 + |p_{n-1}|)$ 
 for the next iteration. 
 In the same way we prove the inequalities for 
 $r_{n-3}, \dots, r_1$ for the next iteration. 
 Finally, since we assumed that  
 $r_1 \leqslant (|q|-1)(1+|p_{n-1}|+ \dots + 
 |p_2|)$ for the current iteration, for the next iteration 
 we have: 
 $|r_0| \leqslant (|q|-1)(1+|p_{n-1}|+|p_{n-2}|+
   \dots + |p_1|)\leqslant (|q|-1)^2$. 
 So, the algorithm requires 
 only a finite amount of memory. The same remains true
 for an algorithm recognizing the addition.  
 In order to get the analogue of Proposition \ref{gcdpqequalone_prop} for $n>2$, one should 
 simply change $p$ to $p_1$. 
 Also, clearly, there is a root $\alpha$ of 
 polynomial $t(x)$ for which 
 $|\alpha| \leqslant \sqrt[n]{|q|}$. 
 
 We call all presentations
 $\psi_{\overline{p},q}$ satisfying the conditions
 $1+|p_{n-1}|+\dots+|p_1|<|q|$ and 
 $\gcd(p_1,q)=1$  
 Nies--Semukhin FA--presentations. 
 So, in exactly the same way as for $n=2$,  
 we obtain that for every Nies--Semukhin FA--presentation 
 the language $L_\xi$ is not regular. 
 Let $\eta_i = [x^i]_\sim$ for $i=1,\dots,n-1$. 
 It is clear that $\mathbb{Z}^n$ is the direct sum 
 of its cyclic subgroups $\langle \xi \rangle$ and 
 $\langle \eta_1 \rangle, \dots, \langle \eta_{n-1} \rangle$. 
 Let $L_{\eta_i} = \psi_{\overline p,q}^{-1} 
 (\langle \eta_i \rangle )$ for $i=1,\dots, n-1$. 
 Similarly to the case $n=2$ we obtain  
 that $L_\xi = \{w \in \mathrm{Dom}_{\overline p,q} \, | \, 
 \varphi_{x} (w) \in  L_{\eta_1} \}$.  
 The inclusion $L_\xi \subseteq 
 \{w \in \mathrm{Dom}_{\overline p,q} \, | \, 
  \varphi_{x} (w) \in L_{\eta_1}\}$ is 
 straightforward. In order to prove the 
 inclusion $\{w \in \mathrm{Dom}_{\overline p,q} \, | \, 
 \varphi_{x}(w) \in L_{\eta_1}\} \subseteq L_\xi$ we 
 notice that if $\psi_{\overline p,q} (w) = 
 [s x^{n-1} + r_{n-2} x^{n-2} + \dots + r_0 ]_\sim$, then  
 $\varphi_x (w) = [x (s x^{n-1} + r_{n-2} x^{n-2} + \dots + r_0)]_\sim= [(r_{n-2} - sp_{n-1})x^{n-1}+ 
 \dots + (r_0 - s p_1)x + sq]_\sim$ which is equal 
 to  $[kx]_\sim$ for some $k \in \mathbb{Z}$ only if 
 $s, r_{n-2},\dots,  r_2 $ 
 and $r_1$ are equal to zero. Thus, the regularity 
 of $L_{\eta_1}$ implies the regularity of 
 $L_\xi$. Therefore, $L_{\eta_1}$ is not regular. 

 Then we consecutively prove that 
 each of the languages 
 $L_{\eta_2}, \dots, L_{\eta_{n-1}}$ is not 
 regular using the observation that 
 $L_{\eta_i} = \{ w \in 
  \mathrm{Dom}_{\overline p, q} \,|\, \varphi_x (w)
   \in L_{\eta_{i+1}}\}$ for $i=1, \dots, n-2$. 
 The inclusion $L_{\eta_i} \subseteq 
 \{w \in \mathrm{Dom}_{\overline p,q} \, | \, 
 \varphi_x(w) \in L_{\eta_{i+1}} \}$ is straightforward. 
 In order to prove the inclusion 
 $L_{\eta_i} \subseteq \{ w \in \mathrm{Dom}_{\overline p, q}
 \, | \, \varphi_x (w) \in L_{\eta_{i+1}}\}$ for 
 $i=1,\dots, n-2$ we notice that  
 if $\psi_{\overline p,q} (w) = 
 [s x^{n-1} + r_{n-2} x^{n-2} + \dots + r_0 ]_\sim$, 
 then $\varphi_x(w) = [kx^{i+1}]_\sim$ for some 
 $k \in \mathbb{Z}$ only if $s=0$ and $r_j =0$ for 
 $j \neq i$.    
 The following theorem 
 generalizes Theorem 
 \ref{niessem_thm1} for the case $n>2$. 
 \begin{theorem}
     For every tuple $\overline{p}= 
     \langle p_1, \dots, p_{n-1} \rangle$ and an integer $q$ 
     for which $1+|p_{n-1}|+\dots+|p_1|<|q|$ 
     the map $\psi_{\overline{p},q} : 
     \mathrm{Dom}_{\overline{p},q} \rightarrow 
     \mathbb{Z}^n$ gives a FA--presentation 
     of $(\mathbb{Z}^n,+)$.	
 	 If $\gcd (p_1,q)=1$, then 
 	 none of the cyclic components of $\mathbb{Z}^n$
 	 and the projections onto these 
 	 components is FA--recognizable with respect to $\psi_{\overline{p},q}$. 
 \end{theorem}
 Clearly, Remarks \ref{mult_by_constant_remark} and  
 \ref{all_cyclic_subgroups_notFA} hold valid 
 also for the case $n>2$. In particular, 
 for every nonzero integer $m$ each cyclic 
 subgroup  
 $\langle  m \eta_1 \rangle, \dots, 
  \langle  m \eta_{n-1} \rangle $ and $\langle m \xi \rangle$
  is not FA--recognizable with respect to 
  a Nies--Semukhin FA--presentation. 
  Furthermore, if a polynomial $t(x)$ is irreducible, 
  then none of the cyclic subgroups of $\mathbb{Z}^n$ 
  is FA--recognizable.    	 
 
\section{FA--recognizable automorphisms of $\mathbb{Z}^n$}
\label{automorphismssection1}

  In this section until the last paragraph we 
  discuss  the case $n=2$.      
  By Proposition 
  \ref{multiplication_by_g_FArec_prop1}, 
  for a polynomial $g \in \mathbb{Z}[x]$, 
  multiplication by $g$ induces a 
  FA--recognizable map 
  $\varphi_g: \mathrm{Dom}_{p,q} 
  \rightarrow \mathrm{Dom}_{p,q}$. 
  Clearly, if $f \sim g$, then
  $\varphi_g = \varphi_f$. Therefore, 
  since every polynomial from 
  $\mathbb{Z}[x]$ is equivalent to a 
  polynomial of degree at most one, we 
  may assume that $g(x)=ax + b$ for 
  $a,b \in \mathbb{Z}$. 
  Let $h(x) = h_1 x + h_2$, for 
  $h_1,h_2 \in \mathbb{Z}$. The 
  equivalence class $[h]_\sim$ is 
  identified with 
  $(h_1,h_2) \in \mathbb{Z}^2$. 
  We have: $g(x)h(x) = (ax+b)(h_1x+h_2)
  = a h_1 x^2 + (a h_2 + b h_1) x + b h_2 
  \sim a h_1 (-px+q) + (a h_2 + b h_1) x + 
  b h_2 = ((b -ap) h_1 + a h_2) x + 
  aq h_1 + b h_2$. 
  Clearly, $\xi=[1]_{\sim}$ and $\eta = [x]_\sim$, 
  already defined in Section 
  \ref{niessemukhinsection}, 
  generate the group $\mathbb{Z}^2$.   
  We denote by $H_1$ and $H_2$ the cyclic 
  subgroups of $\mathbb{Z}^2$ 
  generated by $\eta$ and 
  $\xi$, respectively. 
  Thus, multiplication 
  by $g$ induces an endomorphism of 
  $\mathbb{Z}^2 = H_1 \oplus H_2$ 
  given by a matrix 
  $A = \left(
   \begin{array}{cc} 
     b - ap & a \\
       aq   & b
   \end{array}  
  \right)$.
  The condition that $A \in 
  \mathrm{GL}(2, \mathbb{Z})$ yields 
  the equations 
  $b^2 - ab p - a^2 q = \pm 1$.
  The latter 
  is equivalent  
  to $(2b-ap)^{2}-(p^{2}+4q) a^{2}=\pm4$. 
  Let $c = 2b - ap$. Then we have:
  \begin{equation}  
  \label{matrixAeq1}   
   A = \left(
   \begin{array}{cc} 
     \frac{c-ap}{2} & a \\
       aq   & \frac{c+ap}{2}
   \end{array}  
  \right),
  \end{equation}  
  where $p,q,a$ and $c$ 
  satisfy one of the following two equations:
  \begin{equation}   
  \label{pelleq1}    
    c^2 - (p^2 + 4q)a^2 = \pm 4. 
  \end{equation}  
  For given $p$ and $q$,  
  the trivial solutions of \eqref{matrixAeq1}, 
  $a=0$ and $c=\pm 2$, 
  correspond to the matrices 
  $A=\pm I$.  
  We will assume that 
  $a \neq 0$. 
  Let  $n = p^2 + 4q$.   
  Clearly, nontrivial solutions 
  of \eqref{pelleq1} 
  exist only if $n \geqslant -4$.  
  The following theorem can be verified 
  by direct calculations.   
  \begin{theorem}
  \label{classificationtheorem1}
  For a given $n \geqslant -4$, the   
  matrices $A$ defined  
  by \eqref{matrixAeq1} together with the coefficients 
  $p$ and $q$
  for which $p,q,a$ and $c$ satisfy: 
  $1+|p|<|q|$,  $\gcd(p,q)=1$, $n = p^2 + 4q$, $a \neq 0$ 
  and the equation $c^2 -n a^2 = \pm 4$ 
  are as follows:
    \begin{itemize} 
    \item{For $n=-4$, 
          $A=\pm\left(
          \begin{array}{cc}
          -r & 1\\ 
          -(r^{2}+1) & r
          \end{array}\right)$,
          $p=2r$ and $q=-(r^{2}+1)$, 
          where 
          $r \in (-\infty, -4] \cup 
          [4, +\infty)$ and 
          $r \equiv 0
           \,\,(\mathrm{mod}\,\,2)$.}    
    \item{For $n=-3$, 
         $A=\pm\left(\begin{array}{cc}
         -r & 1\\
         -(r^{2}+r+1) & (r+1)
         \end{array}\right)$
         or 
         $A=\pm\left(\begin{array}{cc}
         -(r+1) & 1\\
         -(r^{2}+r+1) & r
          \end{array}\right)$,
           $p=2r+1$ and $q=-(r^{2}+r+1)$, 
          where $r \in (-\infty, -3] 
          \cup [2, +\infty)$ and
          either $r \equiv 0\,\,(\mathrm{mod}\,\,3)$
          or 
          $r \equiv 2\,\,(\mathrm{mod}\,\,3)$.}    
    \item{For $n=0$, $n=-1$ and 
            $n=-2$, there exist no 
            nontrivial solutions.} 
    \item{For $n=m^2>0$, nontrivial solutions 
    exist only if $n=1$ or $n=4$. 
    For $n=1$,  
    $A=\pm\left(\begin{array}{cc}
    -(2r+1) & 2\\
    -2(r^{2}+r) & (2r+1)
    \end{array}\right)$, $p=2r+1$ and 
    $q=-(r^{2}+r)$, 
    where $r \in (-\infty,-4] \cup 
    [3,+\infty)$. For $n=4$,  
    $A=\pm\left(\begin{array}
    {cc} -r & 1\\
         1-r^{2} & r
    \end{array}\right)$, $p=2r$ and 
    $q=1-r^{2}$, where 
    $r \in (-\infty,-4] \cup 
    [4,+\infty)$ and $r \equiv 0 
    \,\,(\mathrm{mod}\,\,2)$.}    
    \item{For a positive 
    	 nonsquare integer $n$,  
    the equality $n=p^2+4q$ implies 
    that either 
    $n \equiv 0 \,\, (\mathrm{mod}\,\,
    4)$ or $n \equiv 1 \,\, (\mathrm{mod}\,\,
    4)$. For these two cases we have:   
    \begin{itemize}
    \item{For $n =4s$,   
    $A=\pm\left(\begin{array}{cc}
     x-ra & a\\
     a(s-r^{2}) & x+ra
     \end{array}\right) 
    $ or 
    $A=\pm\left(\begin{array}{cc}
     -x-ra & a\\
     a(s-r^{2}) & -x+ra
     \end{array}\right)$,
     $p=2r$ and $q=s-r^{2}$, 
     where $x>0$ and $a>0$ give a  
     solution of Pell's 
     equation or negative Pell's 
     equation: 
     \begin{equation*}     
        x^{2}-sa^{2}=\pm1,   
     \end{equation*}
     and 
     $r$ either satisfies the inequality 
     $|r|<\sqrt{s}-1$ or the 
     inequality 
     $|r|>\sqrt{s+2}+1$. Also, 
     it is required that $\gcd(r,s)=1$ and 
     $r \not\equiv s\,\,(\mathrm{mod}\,\,2)$.}
     \item{For $n\equiv 1\,\, 
     (\mathrm{mod}\,\, 4)$,  
      $A=\pm \left(\begin{array}{cc}
      \frac{c-pa}{2} & a\\
      a\frac{n-p^2}{4} & \frac{c+pa}{2}
      \end{array}\right)$ or 
      $A = \pm \left(\begin{array}{cc}
      \frac{-c-pa}{2} & a\\
      a\frac{n-p^2}{4} & \frac{-c+pa}{2}
      \end{array}\right)     
      $, $p \equiv 1 \,\,
      (\mathrm{mod}\,\,2)$ and 
      $q= \frac{n-p^2}{4}$, where $c>0$ 
      and $a>0$ give a solution of one of the  
      following Pell--type equations: 
      \begin{equation*}
         c^2  - n a^2 = \pm 4,   
      \end{equation*}
      and $p$ either satisfies the 
      inequality $|p|<\sqrt{n}-2$ or
      the inequality $|p|>\sqrt{n+8}+2$. 
      Also, it is required that 
      $\gcd(p,n)=1$.}  
    \end{itemize}    
    }
    \end{itemize}
  \end{theorem}

  \begin{remark} 
     We recall that for a 
     nonsquare integer $n>0$ 
     Pell's equation
     $x^2 - n y^2 = 1$     
     has infinitely 
     many solutions which 
     are recursively 
     generated, using 
     Brahmagupta's identity: 
     $(x_1 ^2 - n y_1 ^2)(x_2 ^2 - n y_2^2) = 
      (x_1 x_2 + n y_1 y_2)^2 - n (x_1 y_2 + y_1 x_2)^2$,      
     from the fundamental solution 
     -- the one for which 
     positive $x$ and $y$ are minimal.      
     The fundamental solution 
     can be found, for example, 
     using continued fraction of 
     $\sqrt{n}$.
     All solutions 
     of negative Pell's equation 
     $x^2 - n y^2 = -1$ 
     are also generated from 
     its fundamental solution. 
     However, solutions of negative 
     Pell's equation do not always 
     exist. The first 54 numbers for 
     which solutions exist 
     are given by 
     the sequence A031396 in OEIS~\cite{oeis}.
     Similarly, for the Pell--type equations 
     $c^2 - n a^2 = 4$ and 
     $c^2 - n a^2 = -4$,  
     all solutions are 
     recursively generated 
     from the fundamental solutions. 
     For the latter equation solutions 
     exist if and only if they exist 
     for the equation $x^2 - n y^2 = -1$.      
     Furthermore, by Cayley's theorem, 
     if the fundamental solution 
     $(u, v)$ of the equation $c^2 - n a^2 = 4$ 
     is odd (i.e., both $u$ and $v$ are odd),
     then 
     $\left((u^2-3)u/2,
     (u^2-1)v/2 \right)$ 
     gives the fundamental solution of
     the equation $x^2 - n y^2 =1$. Similarly, 
     the odd fundamental solution 
     $(u,v)$ of the 
     equation $c^2 - n a^2 = -4$ leads 
     to the fundamental solution 
     $\left((u^2+3)u)/2, 
      ((u^2+1)v)/2\right)$ of the 
      equation $x^2 - n y^2 = - 1$ 
      \cite{pellref1}.  
     If the fundamental solution 
     is even then it is obtained 
     from the fundamental solution 
     of the corresponding Pell's equation 
     by multiplication by $2$.           
  \end{remark}
  \begin{remark} 
  	 We note that Pell's equation already appeared in
  	 the proof that the ring $\left(\mathbb{Z} (\sqrt{n}),\mathbb{Z}, 
  	 +,<,=;\cdot \right)$ has for every positive 
  	 natural number $n$ a semiautomatic presentation
  	 \cite{Semiautomaticstructures17}. 
  	 This is not surprising since the technique 
  	 used in the construction of a such 
  	 semiautomatic presentation 
  	 is similar to the Nies--Semukhin construction
  	 \cite[\S~6]{NiesSemukhin07}.  
  \end{remark}	
   \begin{remark}
   	\label{SpqRemark1}
    For a fixed pair $p$ and $q$, 
    the matrices \eqref{matrixAeq1} with 
    coefficients satisfying \eqref{pelleq1}
    form a submonoid $\mathcal{S}_{p,q}$ in 
    $\mathrm{GL}(2, \mathbb{Z})$. 
    Let $\mathcal{P}$ be the set 
    of all pairs $(p,q)$ for which 
    $1+|p|<|q|$, $\gcd(p,q)=1$ and 
    $n=p^2 + 4q$ is equal to either $-4,-3,1,4$ or 
    a nonsquare positive integer.   
    Then a set of all matrices given by 
    Theorem \ref{classificationtheorem1} together 
    with the matrices $\pm I$ 
    is the union $\mathcal{S} = 
    \bigcup_{(p,q) \in \mathcal{P}} \mathcal{S}_{p,q}$. 
    For different pairs $(p,q), (p',q') \in \mathcal{P}$ 
    we clearly have $\mathcal{S}_{p,q} \cap \mathcal{S}_{p',q'} = \{\pm I \}$. Moreover, it can be verified that each of these submonoids $\mathcal{S}_{p,q}$ is
    isomorphic to one of the  groups:
    $\mathbb{Z}_4$, $ \mathbb{Z}_6$, 
    $\mathbb{Z}_2 \times \mathbb{Z}_2$ and
    $\mathbb{Z} \times \mathbb{Z}_2$. Namely, 
    from Theorem  \ref{classificationtheorem1} 
    we obtain the following. 	
     For $n=-4$, $n=-3$ and 
     $n=1,4$, $\mathcal{S}_{p,q}$  is a finite group isomorphic 
     to $\mathbb{Z}_4$, $ \mathbb{Z}_6$ and 
     $\mathbb{Z}_2 \times \mathbb{Z}_2$, respectively. 
     For a positive nonsquare integer $n$, 
     $\mathcal{S}_{p,q} \cong \mathbb{Z} \times \mathbb{Z}_2$.     
  \end{remark}
  \begin{remark}
  	 Let $(p,q) \in \mathcal{P}$ such 
  	 that the polynomial $t(x)= x^2 + px -q$ is irreducible
  	 in $\mathbb{Z}[x]$.   
  	 One can easily construct an infinite 
  	 family of not FA--recognizable  
  	 automorphisms of $\mathbb{Z}^2$ 
  	 with respect to the representation 
  	 $\psi_{p,q}$. 
  	 Let $A = 
  	 \left(\begin{array}{cc}
  	 a_{11} & a_{12} \\
  	 a_{21} & a_{22}
  	 \end{array}\right)
  	 \in \mathcal{S}_{p,q}$. For a 
  	 matrix $A' = A+D$, where $D=\left(\begin{array}{cc}
  	 k\ell & kn\\
  	 m\ell & mn
  	 \end{array}\right)$ is a nonzero singular matrix, 
  	 $\det A' = \det A$ iff $m(a_{11}n+a_{12}\ell)+k(a_{21}n+a_{22}\ell)=0$. 
  	 The latter equation admits infinitely many solutions for 
  	 $k,l,m$ and $n$. 
  	 Since $A$ is FA--recognizable with respect to 
  	 $\psi_{p,q}$, assuming that $A'$ is FA--recognizable
  	 with respect to $\psi_{p,q}$, 
  	 we  get that $D = A' - A$ must be FA--recognizable with respect to 
  	 $\psi_{p,q}$. But 
  	 $D$ is not FA--recognizable (see Remark 
  	 \ref{degenerateendomorphisms_remark1}), so 
  	 $A'$ is not FA--recognizable.    
  \end{remark} 	
  \begin{remark}
  \label{Heisenberg_remark1}    
     There exist automorphisms of 
     $\mathbb{Z}^2$ which are not FA--recognizable  
     with respect to every representation 
     $\psi_{p,q}$, $(p,q) \in \mathcal{P}$.
     For example, all automorphisms 
     of $\mathbb{Z}^2$ given by the matrices 
     $T_n = \left(\begin{array}{cc}
      1 & 0 \\
     n &  1
     \end{array}\right)$ for nonzero integer $n$ are not 
     FA--recognizable. 
     This follows from 
     the fact that $I$ is FA--recognizable but, 
     by Remark \ref{mult_by_constant_remark},  
     the endomorphisms 
     $T_n - I$
     for $n \neq 0$ are not FA--recognizable.        
     In particular, none of the representations 
     $\psi_{p,q}$, $(p,q) \in \mathcal{P}$ can be used
     to construct 
     a Cayley automatic representation 
     for the Heisenberg group 
     $\mathcal{H}_3(\mathbb{Z}) \cong 
     \mathbb{Z}^2 \rtimes_{T_1} \mathbb{Z}$.     
  \end{remark}	
  \begin{remark}
  	 We note that for two conjugate matrices 
  	 $A$ and $B = TAT^{-1}$ 
  	 in $\mathrm{GL}(2,\mathbb{Z})$ the
  	 groups 
  	 $\mathbb{Z}^2 \rtimes_A \mathbb{Z}$ 
  	 and $\mathbb{Z}^2 \rtimes_B \mathbb{Z}$ 
  	 are isomorphic.     
     An algorithm for solving conjugacy  
     problem in 
     $\mathrm{GL}(2,\mathbb{Z})$ is described in \cite{congGL2Z}; see also 
     an algorithm for solving conjugacy 
     problem in $\mathrm{SL}(2,\mathbb{Z})$ 
     using continued fractions 
     \cite[\S~7.2]{Karpenkovbook}.
     It can be verified that for the cases $n=-4,-3,1,4$ 
     each of the matrices from Theorem \ref{classificationtheorem1}  
     is conjugate to one of the following matrices
     in $\mathrm{GL}(2,\mathbb{Z})$: 
     $\left(\begin{array}{cc}
     0 & -1\\
     1 & 0
     \end{array}\right), \left(\begin{array}{cc}
     1 & 1\\
     -1 & 0
     \end{array}\right), \left(\begin{array}{cc}
     0 & 1\\
     -1 & -1
     \end{array}\right),
     \left(\begin{array}{cc}
     1 & 0\\
     0 & -1
     \end{array}\right)
     $ and 
     $
      \left(\begin{array}{cc}
      0 & 1\\
      1 & 0
      \end{array}\right)
     $.    
     If $n$ is a positive nonsquare integer, 
     every matrix from Theorem \ref{classificationtheorem1},
     which is in $\mathrm{SL}(2,\mathbb{Z})$,  
     is Anosov.       
     Moreover, in this case, for a pair
     $(p,q) \in \mathcal{P}$ satisfying 
     $n = p^2 + 4q$ 
     the matrices from $\mathcal{S}_{p,q}$ 
     generate infinitely many conjugacy classes 
     in $\mathrm{GL}(2,\mathbb{Z})$. The latter 
     immediately follows from the observation 
     that for different values of $c$, 
     which is the trace of the matrix \eqref{matrixAeq1}, 
     we have different conjugacy classes.         
 \end{remark}	
 
  Similarly to the case $n=2$, one can get a family of 
  FA--recognizable automorphisms 
  $A \in \mathrm{GL} (n, \mathbb{Z})$ with respect to the Nies--Semukhin FA--presentations $\psi_{\overline p,q}$ 
  of $\mathbb{Z}^n$ for $n>2$. 
  By Proposition 
  \ref{multiplication_by_g_FArec_prop1} (its 
  analogue clearly holds also for the case $n>2$),
  multiplication by a polynomial $g \in \mathbb{Z}[x]$ 
  induces a FA--recognizable map 
  $\varphi_g : \mathrm{Dom}_{\overline p, q} 
   \rightarrow \mathrm{Dom}_{\overline p, q}$; 
   also, equivalent polynomials $f \sim g$ induce the 
   same map: $\varphi_f = \varphi_g$.  
   So, we may assume that $g(x)= a_{n-1} x^{n-1} + 
   \dots + a_0$ for $a_{n-1}, \dots, a_0 \in \mathbb{Z}$.
   The matrix $A \in \mathrm{GL}(n, \mathbb{Z})$ 
   corresponding to the linear map $\varphi_g$ 
   depends on the tuples 
   $\langle p_1, \dots, p_{n-1} \rangle$, 
   $\langle a_0, \dots, a_{n-1}\rangle$ and 
   the integer $q$. 
   In this paper we do not give full classification 
   of all such matrices for the case $n>2$.


 \section{Conclusion and Open Questions}
 \label{conclusionsec}
 
  In this paper we generalize the Nies--Semukhin 
  FA--presentation of $\left(\mathbb{Z}^2,+\right)$, 
  originally constructed for the polynomial 
  $x^2+ x -3 $, to a polynomial $x^2 + px - q$ 
  such that $1+|p|<|q|$ and $\gcd(p,q)=1$.    
  We also show how this construction is 
  generalized for $\left(\mathbb{Z}^n,+\right), n>2$. 
  Based on this, we construct a new family 
  of Cayley automatic representations 
  of groups $\mathbb{Z}^n \rtimes_A \mathbb{Z}, 
  A \in \mathrm{GL}(n,\mathbb{Z})$ 
  that violate the basic property known for 
  standard representations -- projections 
  $p_i : \mathbb{Z}^n \rightarrow \mathbb{Z}^n, 
  i=1,\dots,n$ are FA--recognizable, 
  i.e., the property b)
  in Section \ref{introsec}.   
  For $n=2$ we describe the set of matrices 
  $\mathcal{S} \subseteq \mathrm{GL}(2,\mathbb{Z})$ 
  corresponding to this family of nonstandard 
  representations and show its connection with
  Pell's equation.  
  Let us pose the following questions that are 
  apparent from the results of this paper. 
  \begin{itemize}
  \item{Is there a nonstandard representation, e.g., 
  	preserving the property a) and violating the property b), 
    for the Heisenberg group $\mathcal{H}_3 (\mathbb{Z})$?}
  \item{What is the set of conjugacy classes of the 
        set of matrices $\mathcal{S}$ in $\mathrm{GL}(2,\mathbb{Z})$?}	
  \item{Is there any Anosov 
  	    $A \in \mathrm{SL} (2,\mathbb{Z})$ 
  	    which is not FA--recognizable with 
  	    respect to every 
  	    Nies--Semukhin FA--presentation 
  	    $\psi_{p,q}, (p,q) \in \mathcal{P}$?  
       } 
  \end{itemize}

\bibliographystyle{splncs04}

\bibliography{lata2020_arxiv}

\end{document}